\documentclass{article}

\usepackage{amscd, amsthm}
\usepackage[footnotesize]{caption}
\usepackage{amsfonts}
\usepackage{amssymb}
\usepackage{amsmath}
\usepackage{latexsym}
\usepackage{amsthm}
\usepackage{eucal}
\usepackage{eufrak}
\usepackage[dvips]{graphicx}
\usepackage{psfrag}
\usepackage{pstricks}

\author{F. Cagliari$^1$,\quad B. Di Fabio$^1$,\quad M. Ferri$^{1,2}$\footnote{Corresponding author. E-mail Address:
\texttt{ferri@dm.unibo.it}}\\\vspace{-0.3cm} \small{$^{1}$
Dipartimento di Matematica, Universit\`a di
Bologna,}\\
\small{P.zza di Porta S. Donato 5, I-$40126$ Bologna,
Italia}\\\vspace{-0.3cm} \small{$^{2}$ARCES, Universit\`a di Bologna,}\\
\small{via Toffano $2/2$, I-$40135$ Bologna, Italia}}

\title{One-Dimensional Reduction of Multidimensional Persistent Homology}

\newcommand{\R}{\mathbb R}
\newcommand{\N}{\mathbb N}
\newcommand{\Z}{\mathbb Z}

\newcommand{\bi}{\vec b}
\newcommand{\li}{\vec l}

\newcommand{\f}{\varphi}
\newcommand{\fr}{\vec f}
\newcommand{\psr}{\vec \psi}

\newcommand{\vhr}{\vec \varphi}

\newcommand{\C}{{\mathcal C}}

\newtheorem{prop}{Proposition}
\newtheorem{definition}{Definition}
\newtheorem{lemma}{Lemma}
\newtheorem{theorem}{Theorem}
\newtheorem{remark}{Remark}
\newtheorem{cor}{Corollary}

\begin{document}
\maketitle

\begin{abstract}
A recent result on size functions is extended to higher homology
modules: the persistent homology based on a multidimensional
measuring function is reduced to a 1-dimensional one. This leads
to a stable distance for multidimensional persistent homology.
Some reflections on $i$-essentiality of homological critical
values conclude the paper.
\end{abstract}

\vspace{1cm}

{\bf Keywords}: Size function, measuring function, rank invariant,
pattern recognition, $i$-essentiality.\\

\section{Introduction}
Topological Persistence started {\it ante litteram} in 1991 with
P. Frosini, who introduced the concept of Size Function
\cite{Fr91},\cite[Sect. 8.4]{KaMiMr04}, a topological-geometrical
tool for describing, analyzing and comparing shapes. This was
actually the origin of rather large experimental research
(\cite{UrVe97,VeUr96,VeUrFrFe93}). Size functions were generalized
by the same School in two directions: Size Homotopy Groups
\cite{FrMu99} (already in a multidimensional setting!) and Size
Functor \cite{CaFePo01}.

At about the same time, Persistent Homology was independently
introduced \cite{EdLeZo00,EdLeZo02} (see also \cite{DeGh,DeRGh}).
All these theories have substantially the same target: shape
recognition. They are constructed on some topological features of
lower level sets of a continuous real-valued function defined on
the object of interest. They also share an important advantage
with respect to other methods of pattern recognition: they capture
qualitative aspects of shape in a formal quantitative way; so,
they turn out to be particularly suited to the analysis of
``natural'' shapes (blood cells, signatures, gestures, melanocytic
lesions, \ldots). Retrospectively, a size function is identifiable
with the rank of a 0-th persistent homology module, while the
first persistent homology module is the Abelianization of the
first size homotopy group \cite{FrMu99}, and the size functor
\cite{CaFePo01} is a functorial formalization of the direct sum
of persistent homology modules.

The results obtained recently, involving the construction of size
functions related to multidimensional measuring functions, lead us
to the same generalization to persistent homology modules, which
is the goal of this paper. As far as Size Theory is concerned,
the main reason for such a generalization is that there are shape
features, that have a multidimensional nature (such as color) and
whose description can be done necessarily by a multidimensional
measuring function. Moreover, there are shapes, which cannot be
discriminated by $n$ size functions related to $n$ different
real-valued measuring functions, but can be distinguished by the
size function related to the $n$-dimensional measuring function of
which those are the components (see Section \ref{examples}).
As mentioned in \cite[Section 2.5]{Gh08}, the study of multidimensional
persistence has strong motivations, but some objective obstacles.
This paper wants to pave a way out of these difficulties.

After recalling some basic notions about multidimensional size
functions and $1$-dimensional persistent homology in Section
\ref{basic}, we adapt the arguments of \cite{BiCeXX} to
multidimensional persistent homology in Section \ref{reduction},
for proving our main result (Theorem \ref{th:main}). This is a
reduction theorem, which takes the detection of discontinuity
points back to the case of $1$-dimensional persistent homology.
This seems to overcome the pessimistic final considerations of
\cite[Section 6]{CaZo06} on the structure of the functions
$\rho_{X,i}$. In fact, although the sets, on which the functions
are constant, are much more complicated than the triangles typical
of the 1-dimensional case, they reduce to them when properly
``sliced'' by a suitable foliation. Stable distances on the leaves
of the foliation define (and approximate) a global distance for rank invariants.
Examples and further remarks on a different kind of reduction conclude the paper.

\section{Basic notions}\label{basic}
In the first part of this section we'll recall briefly the concept
of multidimensional size functions and we'll state the theorem
that gives us the tools to calculate them (Theorem
\ref{szfreduction}). It asserts, indeed, that a suitable planes'
foliation of a $2n$-dimensional real space makes an
$n$-dimensional size function equal to a $1$-dimensional in
correspondence of each plane \cite{BiCeXX}. In the second part we
shall review the definitions of persistent homology module and
related concepts \cite{CoEdHa05}.

\subsection{Multidimensional Size Functions and 1-dimensional reduction}
In Multidimensional Size Theory, any pair $(X, \fr)$, where $X$ is
a non-empty compact and locally connected Hausdorff space, and
$\fr = (f_1, \ldots, f_n): X \rightarrow \R^n$ is a continuous
function, is called a \emph{size pair}. The function $\fr$ is
called an $n$-\emph{dimensional measuring function}. The following
relations $\preceq$ and $\prec$ are defined in $\R^n$: for $\vec u
= (u_1, \ldots, u_n)$ and $\vec v = (v_1, \ldots, v_n)$, we say
$\vec u \preceq \vec v$ (resp. $\vec u \prec \vec v$) if and only
if $u_j \leq v_j$ (resp. $u_j < v_j$) for every index $j = 1,
\ldots, n$. For every $n$-tuple $\vec u = (u_1, \ldots, u_n) \in
\R^n$, let $X\langle\fr \preceq \vec u\rangle$ be the set $\{P \in
X: f_j(P) \leq u_j, j = 1, \ldots, n\}$ and let $\Delta^+$ be the
open set $\{(\vec u, \vec v) \in \R^n \times \R^n: \vec u \prec
\vec v\}$.

\begin{definition}
For every $n$-tuple $\vec v = (v_1, \ldots, v_n) \in \R^n$, we
say that two points $P, Q \in X$ are $\langle\fr \preceq
\vec v\rangle$-\emph{connected} if and only if a connected subset
of $X\langle\fr \preceq \vec v\rangle$ exists, containing $P$ and
$Q$.
\end{definition}

\begin{definition}
The \emph{(n-dimensional) size function} associated with
the size pair $(X, \fr)$ is the function $\ell_{(X, \fr)}: \Delta^+
\rightarrow \N$, defined by setting $\ell_{(X, \fr)}(\vec u, \vec
v)$ equal to the number of equivalence classes in which the set
$X\langle\fr \preceq \vec u\rangle$ is divided by the $\langle\fr
\preceq \vec v\rangle$-connectedness relation.
\end{definition}

An analogous definition for multidimensional persistent homology
will be given in Definition \ref{setdef}.

The main goal of \cite{BiCeXX} for size functions, and of the
present paper for persistent homology, is to reduce computation
from the multidimensional to the 1-dimensional case. This is
possible through particular foliations of $\R^n$ by half-planes.
They are determined by what are called ``admissible'' vector
pairs.

\begin{definition}
For every unit vector $\li = (l_1, \ldots, l_n)$ in $\R^n$ such
that $l_j > 0$ for $j = 1, \ldots, n$, and for every vector $\bi =
(b_1, \ldots, b_n)$ in $\R^n$ such that $\overset{n}{\underset{j =
1}\sum}b_j = 0$, we shall say that the pair $(\li, \bi)$ is
\emph{admissible}. We shall denote the set of all admissible pairs
in $\R^n \times \R^n$ by $Adm_n$. Given an admissible pair $(\li,
\bi)$, we define the half-plane $\pi_{(\li, \bi)}$ in $\R^n \times
\R^n$ by the following parametric equations:
$$
\left\{
\begin{array}{c}
\vec u = s\li + \bi\\
\vec v = t\li + \bi
\end{array}
\right.
$$
for $s, t \in \R$, with $s < t$.
\end{definition}

The motivation for the previous definition is the fact that for
every $(\vec u,\vec v)\in\Delta^+$ there exists exactly one
admissible pair $(\li,\bi)$ such that $(\vec u,\vec
v)\in\pi_{(\li, \bi)}$ \cite[Prop.1]{BiCeXX}. The following Lemma
is substantially contained in the proof of \cite[Thm. 3]{BiCeXX}.

\begin{lemma}\label{setlemma}
Let $(\li, \bi)$ be an admissible pair and $g:
X \rightarrow \R$ be defined by setting
$$g(P) = {\underset{j = 1, \ldots, n}\max}\left\{\frac{f_j(P) -
b_j}{l_j}\right\}.$$ Then, for every $(\vec u, \vec v) = (s\li +
\bi, t\li + \bi) \in \pi_{(\li, \bi)}$, the following equalities
hold:
$$X \langle \fr \preceq \vec u \rangle= {X}\langle g\leq s\rangle,\qquad {X}\langle\fr\preceq\vec v\rangle=
X\langle g\leq t\rangle$$
\end{lemma}
\begin{proof}
For every $\vec u = (u_1, \ldots, u_n) \in \R^n$, with $u_j = sl_j
+ b_j, j = 1, \ldots, n$, it holds that

\begin{eqnarray*}
X \langle \fr \preceq \vec u \rangle& = &\{P \in X: f_j(P) \leq
u_j,
j = 1, \ldots, n\}\nonumber\\
& = & \{P \in X: f_j(P) \leq sl_j + b_j , j = 1, \ldots, n\}\nonumber\\
& = &\{P \in X: \frac{f_j(P) - b_j}{l_j} \leq s, j = 1, \ldots, n\}\nonumber\\
& = & X \langle g \leq s\rangle
\end{eqnarray*}

Analogously, for every $\vec v = (v_1, \ldots, v_n) \in \R^n$,
with $v_j = tl_j + b_j, j = 1, \ldots, n$, it holds that
$X\langle\fr\preceq\vec v\rangle= X\langle g \leq t\rangle$.
\end{proof}

From that, there follows the main theorem of \cite{BiCeXX}:

\begin{theorem}\label{szfreduction}
Let $(\li, \bi)$ and $g$ be defined as in Lemma \ref{setlemma}.
Then the equality
$$\ell_{(X, \fr)}(\vec u, \vec v) = \ell_{(X, g)}(s, t)$$
holds for every \hbox{$(\vec u, \vec v) = (s\li + \bi, t\li + \bi)
\in \pi_{(\li, \bi)}$}.
\end{theorem}

This is indeed the theorem that we are going to extend, in
Section \ref{reduction}, to persistent homology of all degrees. Its importance
resides in the fact that essential discontinuity points
(``cornerpoints'' in the terminology of Size Theory) are the key
to a stable distance between size functions. Unfortunately,
cornerpoints do not form, in general, discrete sets in the
multidimensional case. This theorem makes it possible to find them
``slice by slice'' with the familiar technique of dimension one. A
practical use is for sampling their sets, so getting bounds for a
stable distance between size functions. Our extension will produce the
same opportunity for persistent homology.

\subsection{1-dimensional Persistent Homology}\label{1Dpersistent}
Given a topological space $X$ and an integer $i$, we denote the
$i$-th singular homology module of $X$ over a field $k$ by
$H_i(X)$.

Next we report two definitions of \cite{CoEdHa05}.

\begin{definition}\label{homcrit}
Let $X$ be a topological space and $f$ a real function on $X$. A
\emph{homological critical value} of $f$ is a real number $a$ for
which there exists an integer $i$ such that, for all sufficiently
small $\varepsilon > 0$, the map $H_i (f^{-1}(-\infty, a -
\varepsilon])\rightarrow H_i (f^{-1}(-\infty, a + \varepsilon])$
induced by inclusion is not an isomorphism.
\end{definition}

This is called an $i$-\emph{essential critical value} in the paper
\cite[Def.2.6]{CaFePo01}, dedicated to the \emph{size functor}, a
contemporary and not too different homological generalization of
size functions.

\begin{definition}
A function $f: X \rightarrow \R$ is \emph{tame} if it has a finite
number of homological critical values and the homology modules $H_i
(f^{-1}(-\infty, a])$ are finite-dimensional for all $i \in \Z$
and $a \in \R$.
\end{definition}

The reader should be warned that there exist other, different
meanings of ``tame'' in the current topological literature.
Actually, ``homologically tame'' might be a better designation for
such a type of function, but we adhere to this already current
definition.

We write $F_i^{u} = H_i(f^{-1}(-\infty, u])$, for all $i \in \Z$,
and for $u < v$, we let $f_i^{u,v}: F_i^{u} \rightarrow F_i^{v}$
be the map induced by inclusion of the lower level set of $u$ in
that of $v$, for a fixed integer $i$. Moreover, we indicate with
$F_i^{u,v} = \mbox{Im} f_i^{u,v}$ the image of $F_i^{u}$ in
$F_i^{v}$, that is called $i$-\emph{th persistent homology module}.

\section{Homological 1-dimensional reduction}\label{reduction}

In this section we define the $i$-th persistent homology module
related to a continuous $n$-dimensional real function
(substantially as in \cite{CaZo06}). Then we show that the sets of
points of $\R^{2n}$, where the modules change, can be obtained by
computing the discontinuity points of persistent homology of a
1-dimensional function defined on particular half-planes which
foliate the $2n$-space.

The first issue arises when one tries to compute the maximum
between the components of a $n$-dimensional real function. In
fact:

\begin{remark}
The maximum of two tame functions is not necessarily a tame
function.
\end{remark}

(We recall that ``tame'' has the meaning defined in Section
\ref{1Dpersistent}.)

As an example, let $f_1, f_2: \R^2 \rightarrow \R$ be two tame
functions defined as
$$
\begin{array}{lr}
f_1(u, v) = \left\{
\begin{array}{ll}
 v - u^2 \sin (\frac{1}{u})& u \neq 0\\
 v & u = 0
\end{array}
\right., & f_2(u, v) = \left\{
\begin{array}{ll}
-v - u^2 \sin (\frac{1}{u})& u \neq 0\\
-v & u = 0
\end{array}
\right.
\end{array}
$$
 and consider the function $$f = \max
(f_1,f_2).$$

\begin{figure}[htbp]
\centering
\begin{minipage}[c]{.40\textwidth}
\centering\setlength{\captionmargin}{20pt}%
\psframebox[framesep=0pt,
linecolor=gray]{\includegraphics[scale=0.14, angle=0]{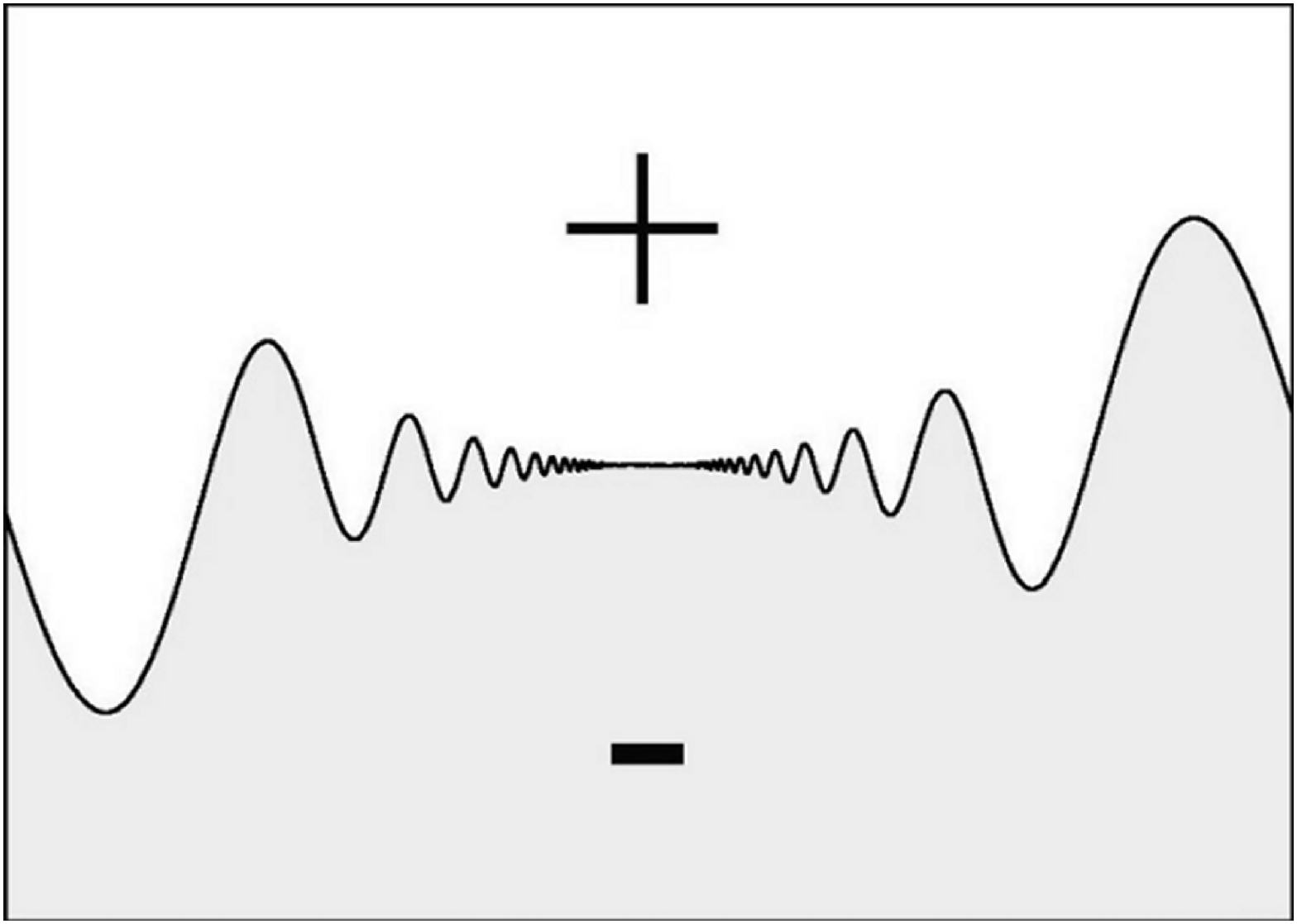}}
\caption{\footnotesize{Lower level set of $f_1$ (grey area - one
connected component).}}
\end{minipage}%
\begin{minipage}[c]{.40\textwidth}
\centering\setlength{\captionmargin}{20pt}%
\psframebox[framesep=0pt,
linecolor=gray]{\includegraphics[scale=0.14, angle=0]{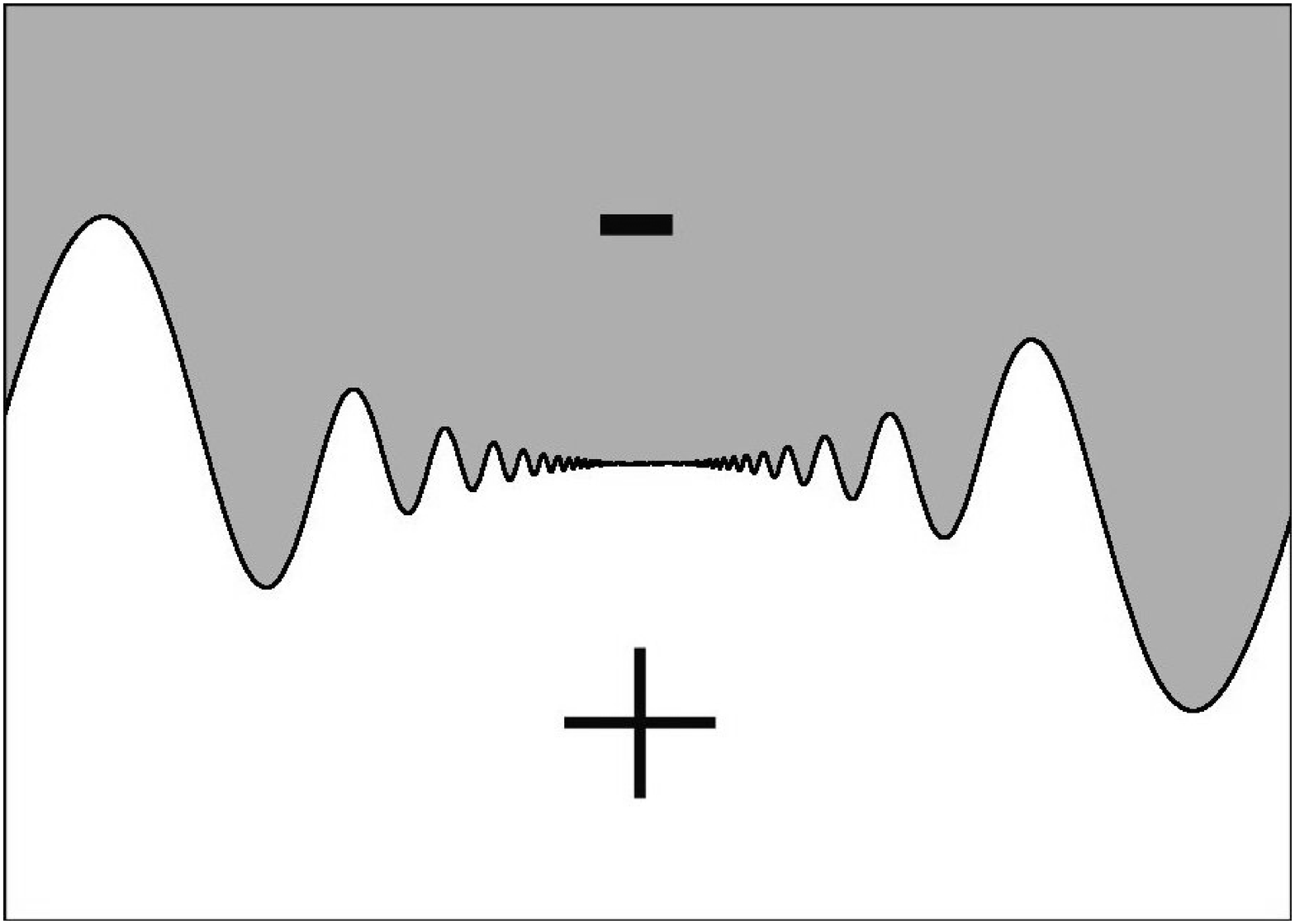}}
\caption{\footnotesize{Lower level set of $f_2$ (grey area - one
connected component).}}
\end{minipage}%
\vspace{3mm}
\begin{center}
\psframebox[framesep=0pt,
linecolor=gray]{\includegraphics[scale=0.14, angle=0]{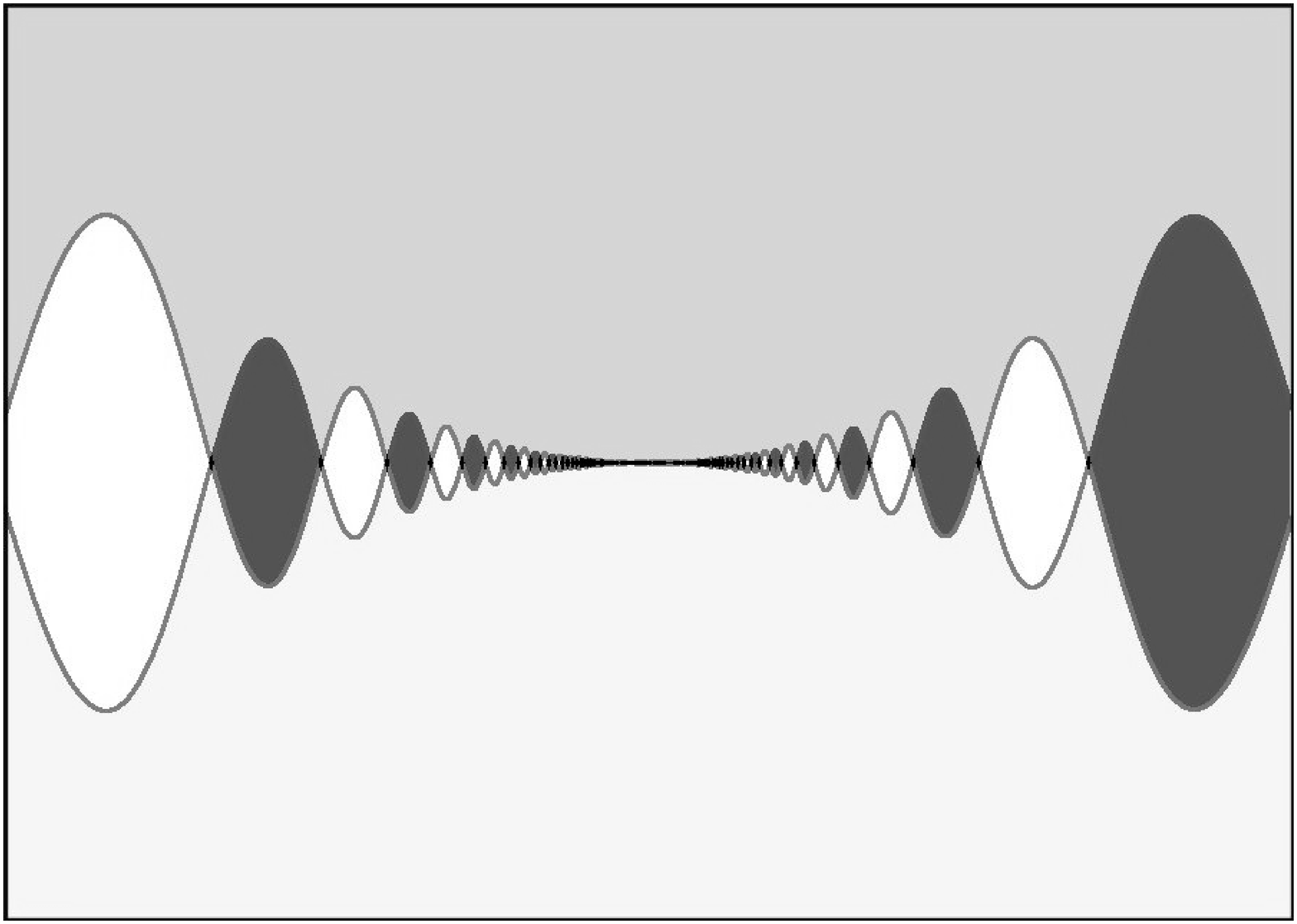}}
\caption{\footnotesize{Lower level set of $f$ (dark zone -
infinitely many connected components).}}
\end{center}
\end{figure}\label{notame}

Then, as we can see in Figure \ref{notame}, $f$ is not tame, since
$H_0(f^{-1}(-\infty, 0])$ is an infinite-dimensional module.

Given this fault related to tame functions, the solution we
propose is to introduce the following concept.

\begin{definition}
 Let $X$ be a topological space and $\fr: X \rightarrow
\R^n$ a continuous function on $X$. We shall say that $\fr$ is
\emph{max-tame} if, for every admissible pair $(\li, \bi)$, the
function $g(P) = {\underset{j = 1, \ldots, n}\max}\left\{\frac{f_j(P)
- b_j}{l_j}\right\}$ is tame.
\end{definition}
Choosing a measuring function on $X$ as above, let us define the
multidimensional persistent modules.
\begin{definition}\label{setdef}
Let $\fr: X \rightarrow \R^n$ be a max-tame function. For each
homology degree $i \in \Z$ we put $F_i^{\vec u} =
H_i({\fr}^{-1}(\overset{n}{\underset{j = 1}\prod}(-\infty,
u_j]))$, for all $\vec u \in \R^n$. For $\vec u \preceq \vec v$ we
let $f_i^{\vec u, \vec v}: F_i^{\vec u} \rightarrow F_i^{\vec v}$
be the map induced by inclusion of the lower level set of $\vec u$
in that of $\vec v$, for a fixed integer $i$, and call $F_i^{\vec
u, \vec v} = \emph{\mbox{Im}}\ f_i^{\vec u, \vec v}$ the
$i$-\emph{th multidimensional persistent homology module}.
\end{definition}

Note that the rank of $F_i^{\vec u, \vec v}$ is what is
called $\rho_{X, i}(\vec u, \vec v)$ in \cite[Def. 12]{CaZo06}.

Let $g(P) = {\underset{j = 1, \ldots, n}\max}\left\{\frac{f_j(P) -
b_j}{l_j}\right\}$ for a fixed $(\li, \bi) \in Adm_n$, $G_i^{s} =
H_i(g^{-1}(-\infty, s])$, for all $s \in \R$ and $i \in \Z$. For
$s < t$, we let $g_i^{s,t}: G_i^{s} \rightarrow G_i^{t}$ be the
map induced by inclusion of the lower level set of $s$ in that of
$t$, for a fixed integer $i$, and denote $G_i^{s, t} = \mbox{Im}\
g_i^{s, t}$ the $i$-th persistent homology module.

Now we can state and prove the theorem which, in analogy with the
main result of \cite{BiCeXX}, enables us to reduce the computation
of multidimensional persistent homology to the 1-dimensional one.
This is important, not so much for finding the homology modules
themselves point by point, but much more for finding points of
change of the modules.

\begin{theorem}\label{th:main}
Let $(\li, \bi)$ be an admissible pair and $\fr = (f_1, \ldots,
f_n): X \rightarrow \R^n$ a max-tame function. Then, for every
$(\vec u, \vec v) = (s\li + \bi, t\li + \bi) \in \pi_{(\li,
\bi)}$, the following equality
$$F_i^{\vec u,\vec v} = G_i^{s,t}$$
holds for all $i \in \Z$ and $s, t \in \R$ with $s < t$.
\end{theorem}
\begin{proof}
By Lemma \ref{setlemma}, we know that, for every $\vec u = (u_1,
\ldots, u_n) \in \R^n$, with $u_j = sl_j + b_j, j = 1, \ldots, n$,
it holds that
$$\{P \in X, f_j(P) \leq u_j, j = 1, \ldots, n\} = \{P \in X, g(P)
\leq s\}$$ hence
$$\{P \in X, P \in {f_j}^{-1}(-\infty, u_j], j = 1, \ldots, n\} = \{P
\in X, P \in g^{-1}(-\infty, s]\}.$$ It follows that
$$\overset{n}{\underset{j = 1}\bigcap}{f_j}^{-1}(-\infty, u_j] =
g^{-1}(-\infty, s]$$ implying
$$H_i\left(\overset{n}{\underset{j =
1}\bigcap}{f_j}^{-1}(-\infty, u_j]\right) = H_i(g^{-1}(-\infty,
s])$$ for all $i \in \Z$.

Analogously, for every $\vec v = (v_1, \ldots, v_n) \in \R^n$,
with $v_j = tl_j + b_j, j = 1, \ldots, n$, it holds that
$H_i\left(\overset{n}{\underset{j = 1}\bigcap}{f_j}^{-1}(-\infty,
v_j]\right) = H_i(g^{-1}(-\infty, t])$, for all $i \in \Z$. So,
since $f_i^{\vec u,\vec v}$ and $g_i^{s,t}$ have the same domain
and codomain and they are the maps induced by inclusion, we can
conclude that $f_i^{\vec u,\vec v} = g_i^{s,t}$, and the claim
follows.
\end{proof}

\section{Multidimensional matching distance}\label{distance}

According to \cite[Def. 12]{CaZo06}, for a given measuring
function $\fr': X \rightarrow \R^n$, for each homology degree
$i\in\Z$ the rank invariant $\rho'_{X,i}:\Delta^+\to \N$ is
defined as $\rho'_{X,i}(\vec u,\vec v)=rank(F_i^{\vec u, \vec
v})$.

Let $(X,\fr')$, $(Y,\fr'')$ be two size pairs, where $\fr': X
\rightarrow \R^n , \ \fr'': Y \to \R^n$ are max-tame measuring
functions, and $\overline{\rho}'_{X,i},\ \overline{\rho}''_{Y,i}$
be  the respective rank invariants. Let an admissible pair $(\li,
\bi)$ be fixed, and let $g': X \rightarrow \R,\ \  g'':Y \to \R$
be defined by setting
$$g'(P) = {\underset{j = 1, \ldots, n}\max}\left\{\frac{f'_j(P) -
b_j}{l_j}\right\} \ \ \ \ \ \ g''(P) = {\underset{j = 1, \ldots, n}\max}\left\{\frac{f''_j(P) -
b_j}{l_j}\right\}$$

It is well-known for 1-dimensional measuring functions
\cite{LaFr97,FrLa01,CoEdHa05} that the relevant information on the
rank invariants $\rho'_{X,i},\ \rho''_{Y,i}$ of $g'$ and $g''$
respectively is contained, for each degree $i$, in their multisets
of cornerpoints, which are called ``persistence diagrams''. These
are sets of points of the extended plane with multiplicities,
augmented by adding a countable infinity of points of the diagonal
$y=x$: let them be called respectively $C'$ and $C''$. Each
cornerpoint is determined by its coordinates $x<y\le\infty$. The
distance of two cornerpoints is

$$\delta\left((a,b),(c,d)\right)=\min\left\{\max\{|a-c|,|b-d|\},
\max\left\{\frac{b-a}{2}, \frac{d-c}{2}\right\}\right\}$$

It has been proved in \cite{CoEdHa05} that the \emph{matching} (or {\emph{bottleneck}) distance

$$d(\rho'_{X,i}, \rho''_{Y,i})={\underset{\sigma}\min}\ {\underset{P\in C'}\max}\  \delta(P,\sigma(P))$$

\noindent where $\sigma$ varies among all bijections from $C'$ to
$C''$, is stable. Mimicking \cite{BiCeXX} (and recalling that
$\rho'_{X,i}, \rho''_{Y,i}$ vary with $(\li,\bi)$) we can use $d$
to define distances between the rank invariants of the original
multidimensional persistent homologies.

\begin{definition}\label{multid}
Let $(X,\fr')$, $(Y,\fr'')$ be two size pairs and
$\overline{\rho}'_{X,i},\ \overline{\rho}''_{Y,i}$ be  the
respective rank invariants. Then the $i$-th \emph{multidimensional
matching distance} between rank invariants is defined as the
extended distance
$$D(\overline{\rho}'_{X,i}, \overline{\rho}''_{Y,i})=
{\underset{(\li,\bi)\in Adm_n}\sup}\ {\underset{j=1,\ldots,n}\min}\ l_j \cdot \ d(\rho'_{X,i}, \rho''_{Y,i})$$
\end{definition}

Note that $D$ is by construction a global distance, i.e. not
depending on $(\li,\bi)$, but since the coefficients $l_j$ are
$\le 1$, there might be distances $d$, for particular admissible
pairs, which take greater values. An easy corollary of our Theorem \ref{th:main} is the
following, which is the higher degree version of \cite[Cor.
1]{BiCeXX}.

\begin{cor}
For each $i\in\Z$ the identity $\overline{\rho}'_{X,i} \equiv \overline{\rho}''_{Y,i}$
holds if and only if $d(\rho'_{X,i}, \rho''_{Y,i})=0$ for every admissible pair $(\li,\bi)$.
\end{cor}

With the same argument of the analogous Proposition 4 of
\cite{BiCeXX}, it is easy to prove the following inequality
between the multidimensional matching distance and the
1-dimensional one obtained by considering the components of the
measuring functions. That this inequality can be strict, is shown
in Section \ref{examples}.

\begin{prop}\label{components}
Let $(X,\fr),(Y,\vec h)$ be size two pairs with
$\fr=(f_1,\ldots,f_n),\ \vec h=(h_1,\ldots,h_n)$ max-tame
measuring functions. For each $i\in\Z$ and for each $j=1,\ldots,n$
let $\rho^{f_j}_{X,i}, \rho^{h_j}_{Y,i}$ be the $i$-th rank
invariants relative to the components $f_j,\ h_j$ respectively;
let then $\overline{\rho}'_{X,i},\ \overline{\rho}''_{Y,i}$ be the
rank invariants relative to $\fr,\ \vec h$ respectively. Then it
holds that
$$d(\rho^{f_j}_{X,i},\rho^{h_j}_{Y,i})\ \le\ D(\overline{\rho}'_{X,i}, \overline{\rho}''_{Y,i})$$
\end{prop}

The matching distance is known to be stable with respect to
perturbation of 1-dimensional measuring functions \cite[Section
3.1]{CoEdHa05} \cite[Thm. 25]{DAFrLa}. In the multidimensional
setting, the stability of $d$ with respect to an admissible pair
is stated in the following proposition, whose proof is again a
copy of that of \cite[Prop. 2]{BiCeXX}. Here $g',g'':X\to\R$ are
defined in correspondence to $(\li,\bi)$ as at the beginning of
this Section.

\begin{prop}\label{stab}
If $(X,\fr'),\ (X, \fr'')$ are size pairs, with max-tame functions
$\fr', \fr'':X\to\R^n$, and $\max_{P\in X} \Vert
\fr'(P)-\fr''(P)\Vert_\infty \le \epsilon$, then for every
admissible pair $(\li,\bi)$ and for each $i\in\Z$ it holds that
$$d(\rho'_{X,i}, \rho''_{X,i})\le \frac{\epsilon}{\min_{j=1,\ldots,n}\ l_j}$$
with $\li=(l_1,\ldots,l_n)$ and where $\rho'_{X,i},\ \rho''_{X,i}$ are the rank
invariants at degree $i$ of $(X,g'),\ (X,g'')$ respectively.
\end{prop}

By the definition of $D$, every 1-dimensional matching distance
obtained in correspondence of an admissible pair yields a lower
bound for the multidimensional matching distance $D$; a
sufficiently fine sampling by admissible pairs produces
approximations of arbitrary precision of it.

Of course, $\underset{i\in\Z}\max D(\overline{\rho}'_{X,i},
\overline{\rho}''_{Y,i})$ is still a meaningful distance related
to the size pairs $(X,\fr'), (Y,\fr'')$. We plan to study its
relation with the $n$-dimensional natural pseudodistance of
\cite{FrMu99,BiCeXX}.

\section{Examples and Remarks}\label{examples}

We now describe a simple example, which shows that persistent
homology, with respect to a multidimensional measuring function,
is actually stronger than the simple collection of the persistent
homologies with respect to its 1-dimensional components.

In $\R^3$ consider the set $\Omega = [-1, 1]\times[-1,
1]\times[-1, 1]$ and the sphere $\mathcal S$ of equation $u^2 +
v^2 + w^2 = 1$. Let also $\fr = (f_1, f_2): \R^3 \rightarrow \R^2$
be a continuous function, defined as $\fr(u, v, w) = (|u|, |v|).$
In this setting, consider the size pairs $(\C, \vhr)$ and
$(\mathcal S, \psr)$, where $\C =
\partial\Omega$ and $\vhr$ and $\psr$
are respectively the restrictions of $\fr$ to $\C$ and $\mathcal
S$.

In order to compare the persistent homology modules of $\C$ and
$\mathcal S$ defined by $\fr$, we are interested in studying the
half-planes' foliation of $\R^4$, where $\li = (\cos\theta,
\sin\theta)$ with $\theta \in (0, \frac{\pi}{2})$, and $\bi = (a,
-a)$ with $a \in \R$. Any such half-plane is parameterized as
$$
\left\{
\begin{array}{c}
u_1 = s\cos\theta + a\\
u_2 = s\sin\theta - a\\
v_1 = t\cos\theta + a\\
v_2 = t\sin\theta - a
\end{array}
\right.
$$
with $s, t \in \R, s < t$.

In the following, we shall always assume $0 \leq s < t$.

For example, by choosing $\theta = \frac{\pi}{4}$ and $a = 0$,
i.e. $\li = (\frac{\sqrt{2}}{2}, \frac{\sqrt{2}}{2})$ and $\bi =
(0, 0)$, we obtain that
$$g' = \sqrt{2}\max\{\f_1, \f_2\} =
\sqrt{2}\max\{|u|, |v|\},$$ $$g'' = \sqrt{2}\max\{\psi_1, \psi_2\}
= \sqrt{2}\max\{|u|, |v|\}.$$

\begin{figure}[htbp]
\centering
\begin{minipage}[c]{.40\textwidth}
\centering\setlength{\captionmargin}{0pt}%
\includegraphics[scale=0.14, angle=0]{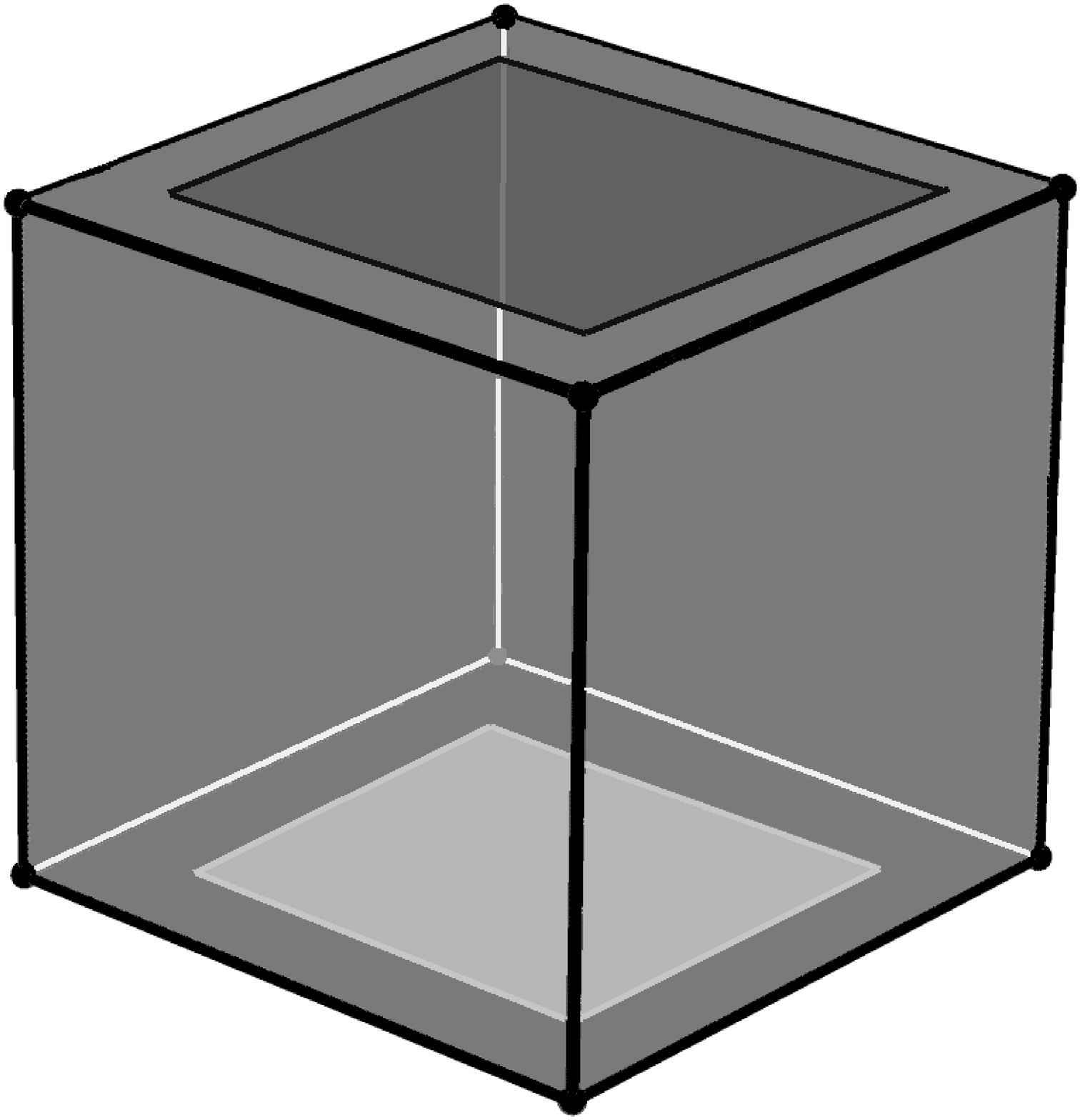}
\end{minipage}%
\hspace{4mm}%
\begin{minipage}[c]{.40\textwidth}
\centering\setlength{\captionmargin}{0pt}%
\includegraphics[scale=0.14, angle=0]{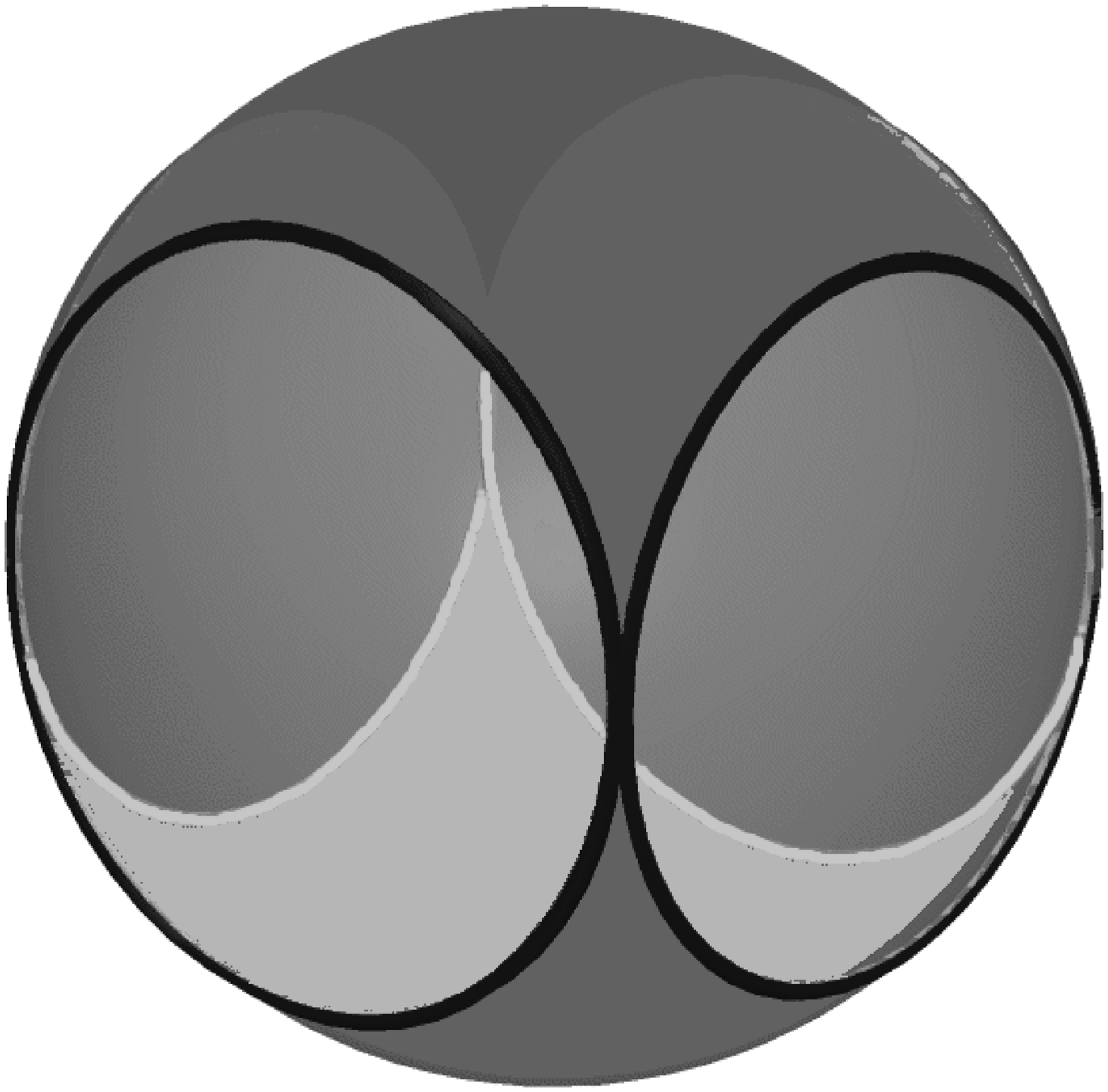}
\end{minipage}
\caption{\footnotesize{Lower level sets $g' \leq 1$ and $g'' \leq
1$.}}
\end{figure}

Let $\rho'_{\C,i}\ \rho''_{\mathcal S,i}$ be the rank invariants of the respective persistent
homologies for $i\in\Z$.
So, writing $G_i^{s}(\C) =
H_i\left(\left(g'\right)^{-1}(-\infty, s]\right)$,
$G_i^{s}(\mathcal S) = H_i\left(\left(g''\right)^{-1}(-\infty,
s]\right)$ and $G_i^{s,t}$ defined as above, we obtain that

\begin{tabular}{lrl}
\!\!\!\!\!\!\!\!\!\!\!\!\!\!\!\!\!\!\!\!\!\!\!\!\!\!\!\!\!\!\!\!\!\!\!\!\!\!\!\!\!\!\!\!\!\!\!\!&
$\left.
\begin{array}{c}
\begin{tabular}{lr}
$G_0^{s,t} (\C) = \left\{
\begin{array}{ll}
0, & s,t<0\\
k^2,& 0\le s<t < \sqrt{2}\\
k, & \mbox{otherwise}
\end{array}
\right.$&
\begin{tabular}{c}
\psfrag{0}{$\scriptstyle{0}$} \psfrag{1}{$\scriptstyle{1}$}
\psfrag{2}{$\scriptstyle{2}$} \psfrag{s}{$\scriptstyle{s}$}
\psfrag{t}{$\scriptstyle{t}$} \psfrag{t=a}{$\scriptstyle{t =
\sqrt{2}}$}
\includegraphics[height=2cm]{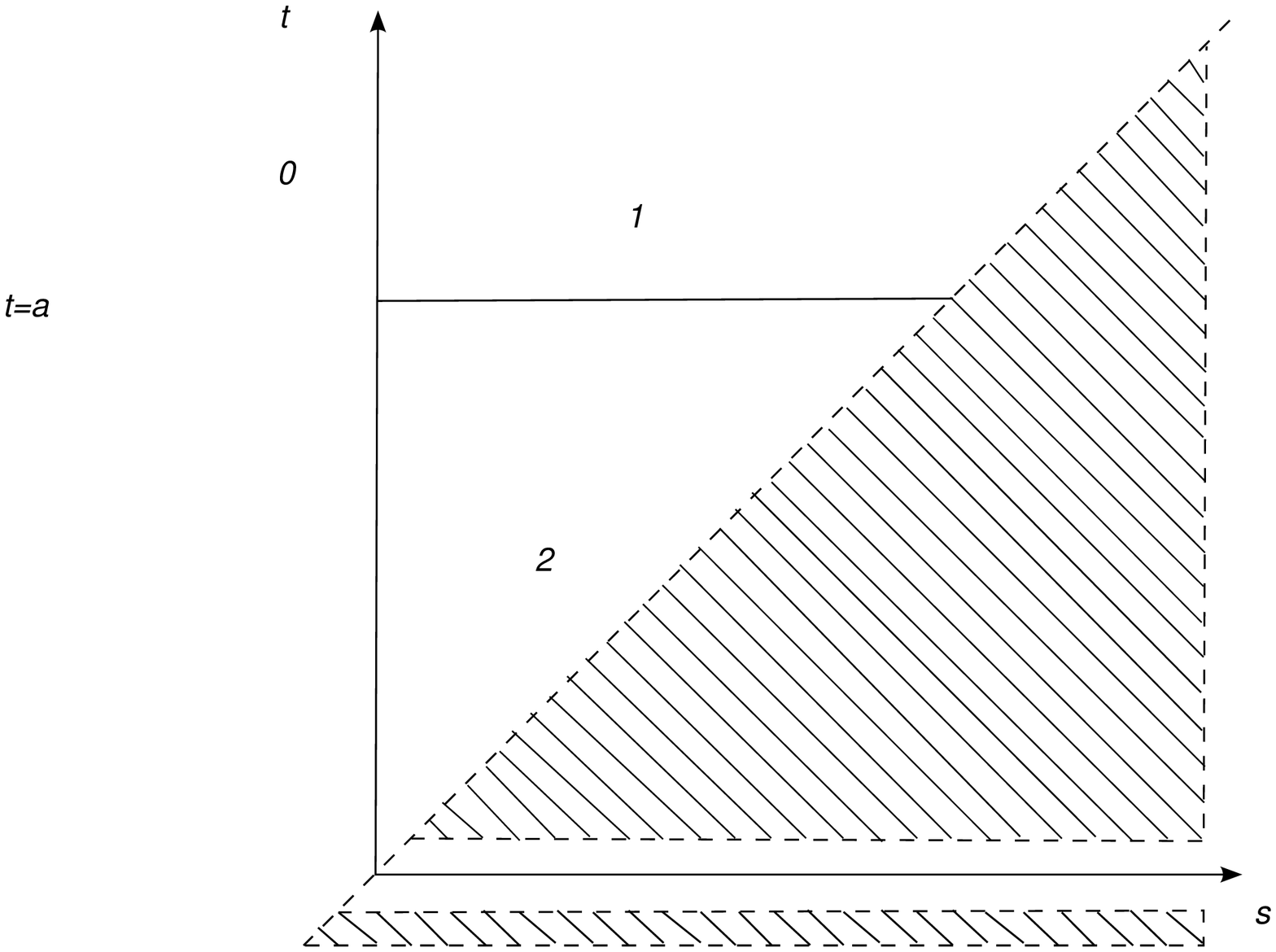}
\end{tabular}
\\
$G_0^{s,t} (\mathcal S) = \left\{
\begin{array}{ll}
0, & s,t<0\\
k^2,&  0\le s<t < 1\\
k, & \mbox{otherwise}
\end{array}
\right.$&
\begin{tabular}{c}
\psfrag{0}{$\scriptstyle{0}$} \psfrag{1}{$\scriptstyle{1}$}
\psfrag{2}{$\scriptstyle{2}$} \psfrag{s}{$\scriptstyle{s}$}
\psfrag{t}{$\scriptstyle{t}$} \psfrag{t=a}{$\scriptstyle{t = 1}$}
\includegraphics[height=2cm]{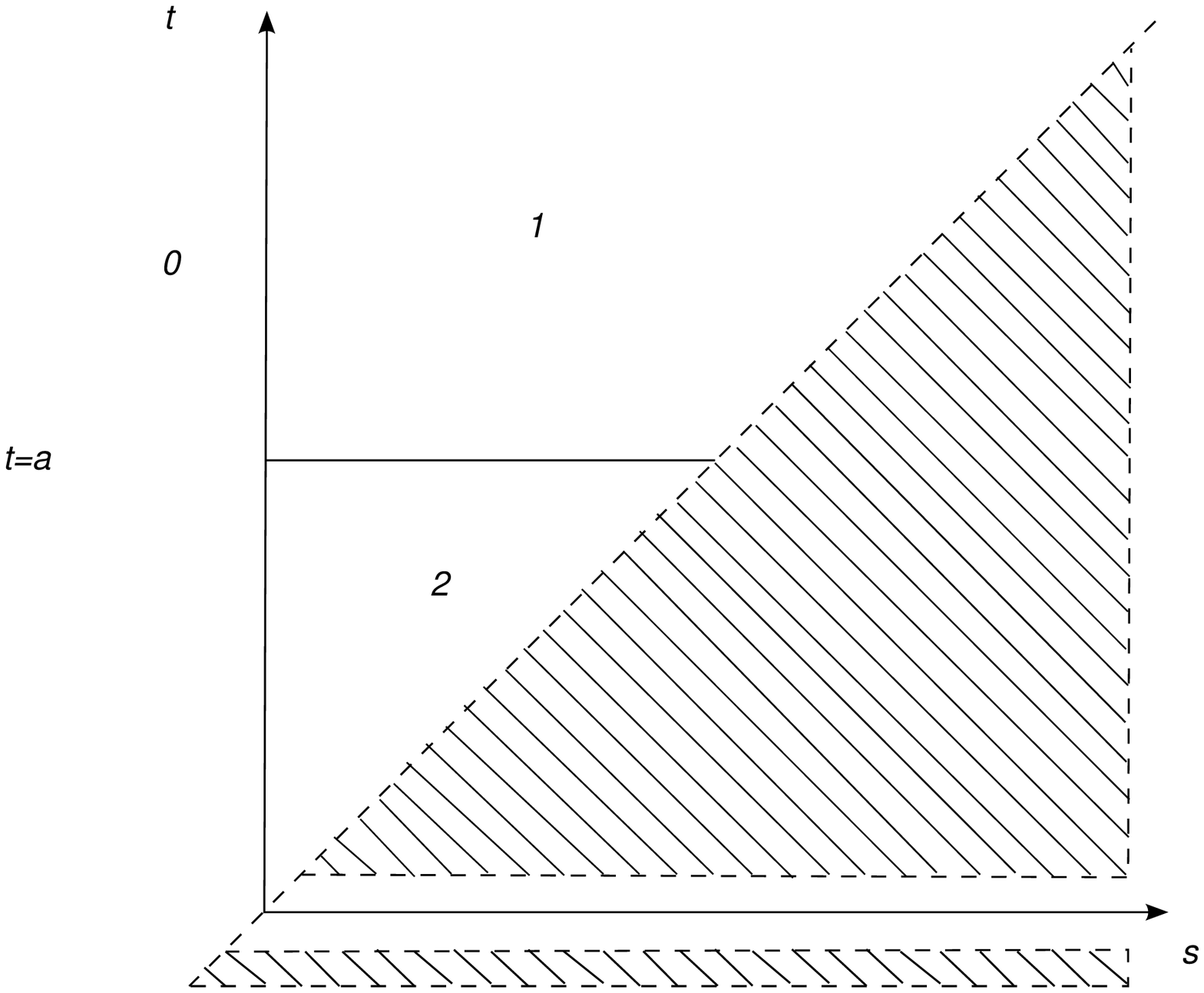}
\end{tabular}\\
\end{tabular}
\end{array}\right\}$&$\Rightarrow D(\rho_{\C,0},\rho_{\mathcal S,0})\geq
 \frac{\sqrt{2}}{2} d(\rho'_{\C,0}\ \rho''_{\mathcal S,0}) = \frac{\sqrt{2}}{2}(\sqrt{2} - 1)$\\
\!\!\!\!\!\!\!\!\!\!\!\!\!\!\!\!\!\!\!\!\!\!\!\!\!\!\!\!\!\!\!\!\!\!\!\!\!\!\!\!\!\!\!\!\!\!\!\!&
$\left.
\begin{array}{c}
\begin{tabular}{lr}
$G_1^{s,t} (\C) = \left.
\begin{array}{lr}
0,& \mbox{for all}\, s, t \in \R
\end{array}
\right.$&
\begin{tabular}{c}
\psfrag{0}{$\scriptstyle{0}$} \psfrag{1}{$\scriptstyle{1}$}
\psfrag{2}{$\scriptstyle{2}$} \psfrag{s}{$\scriptstyle{s}$}
\psfrag{t}{$\scriptstyle{t}$} \psfrag{t=a}{$\scriptstyle{t = 1}$}
\includegraphics[height=2cm]{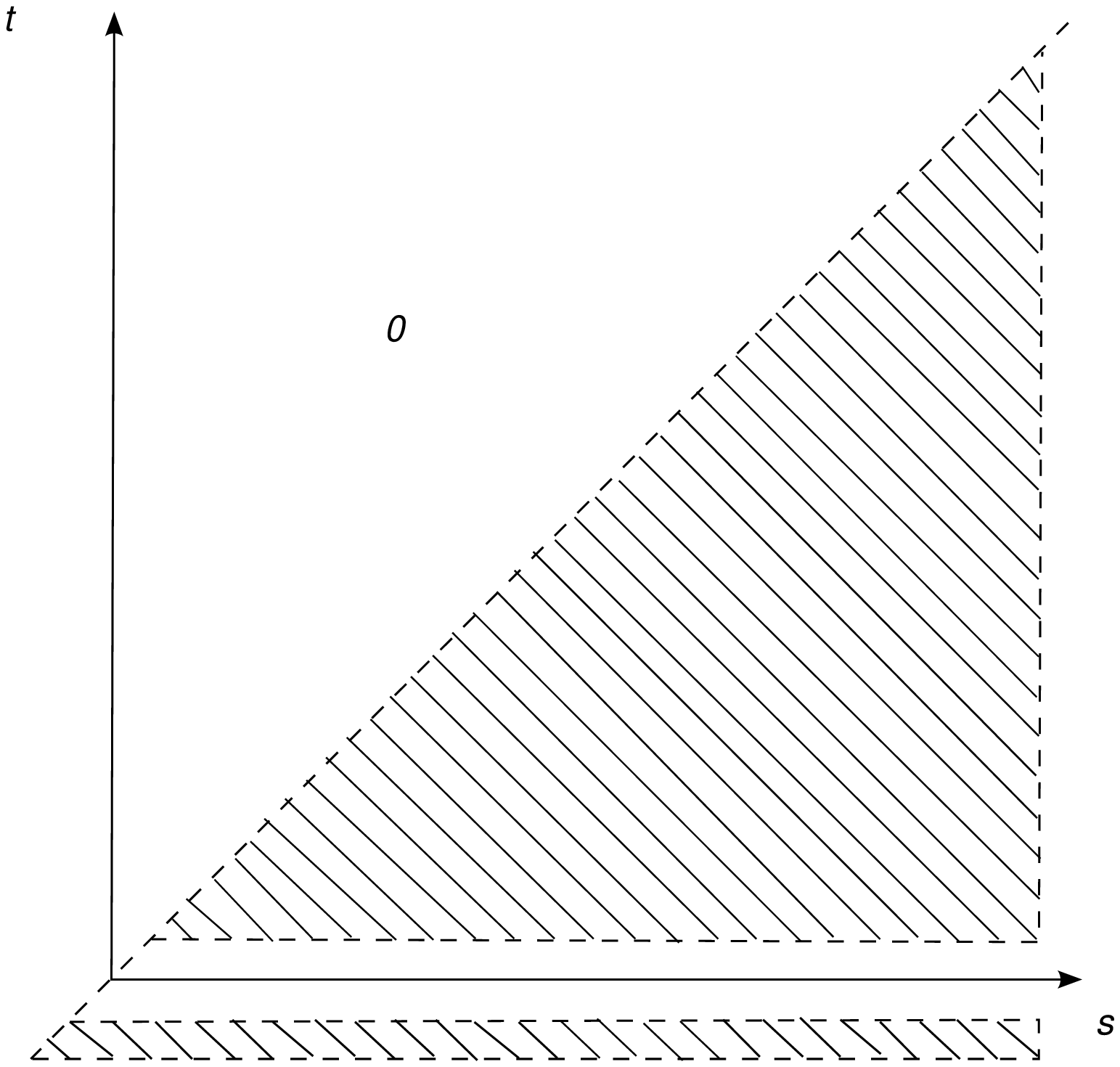}
\end{tabular}\\
$G_1^{s,t} (\mathcal S) = \left\{
\begin{array}{ll}
k^3,& 1 \leq s < t < \sqrt{2}\\
0, & \mbox{otherwise}
\end{array}
\right.$\qquad&
\begin{tabular}{c}
\psfrag{0}{$\scriptstyle{0}$} \psfrag{1}{$\scriptstyle{1}$}
\psfrag{3}{$\scriptstyle{3}$} \psfrag{s}{$\scriptstyle{s}$}
\psfrag{t}{$\scriptstyle{t}$} \psfrag{t=a}{$\scriptstyle{t =
\sqrt{2}}$}
\includegraphics[height=2cm]{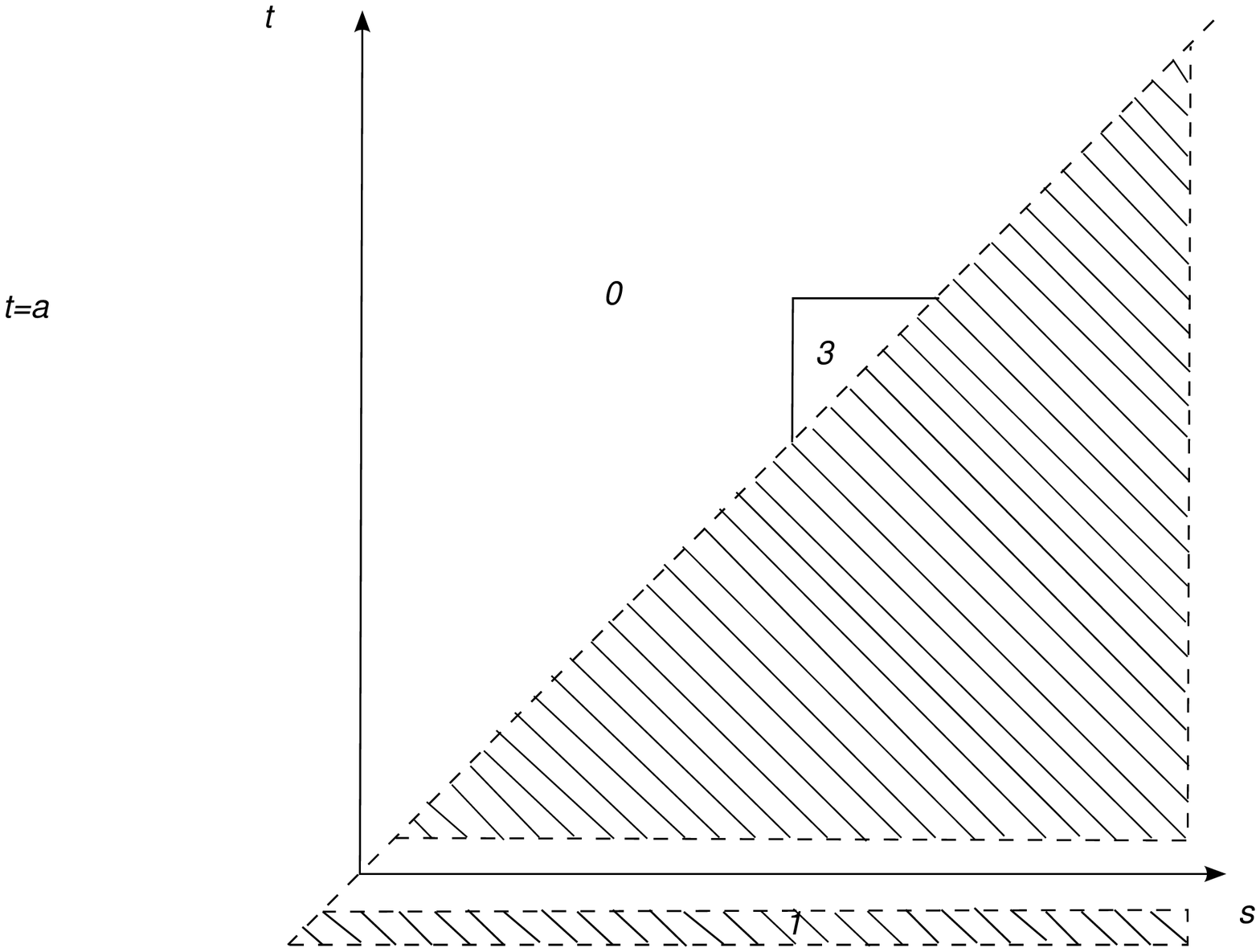}
\end{tabular}\\
\end{tabular}
\end{array}\right\}$ &$\Rightarrow D(\rho_{\C,1},\rho_{\mathcal S,1})\geq
\frac{\sqrt{2}}{2} d(\rho'_{\C,1}\ \rho''_{\mathcal S,1}) = \frac{\sqrt{2}}{2}\left(\frac{\sqrt{2} - 1}{2}\right)$\\
\!\!\!\!\!\!\!\!\!\!\!\!\!\!\!\!\!\!\!\!\!\!\!\!\!\!\!\!\!\!\!\!\!\!\!\!\!\!\!\!\!\!\!\!\!\!\!\!&
$\left.
\begin{array}{c}
\begin{tabular}{lr}
$G_2^{s,t} (\C) = \left\{
\begin{array}{ll}
k,&  \sqrt{2}\le s<t\\
0, & \mbox{otherwise}
\end{array}
\right.$\qquad\,\qquad\,&
\begin{tabular}{c}
\psfrag{0}{$\scriptstyle{0}$} \psfrag{1}{$\scriptstyle{1}$}
\psfrag{2}{$\scriptstyle{2}$} \psfrag{s}{$\scriptstyle{s}$}
\psfrag{t}{$\scriptstyle{t}$}
\psfrag{t=a}{$\scriptstyle{\sqrt{2}}$}
\includegraphics[height=2cm]{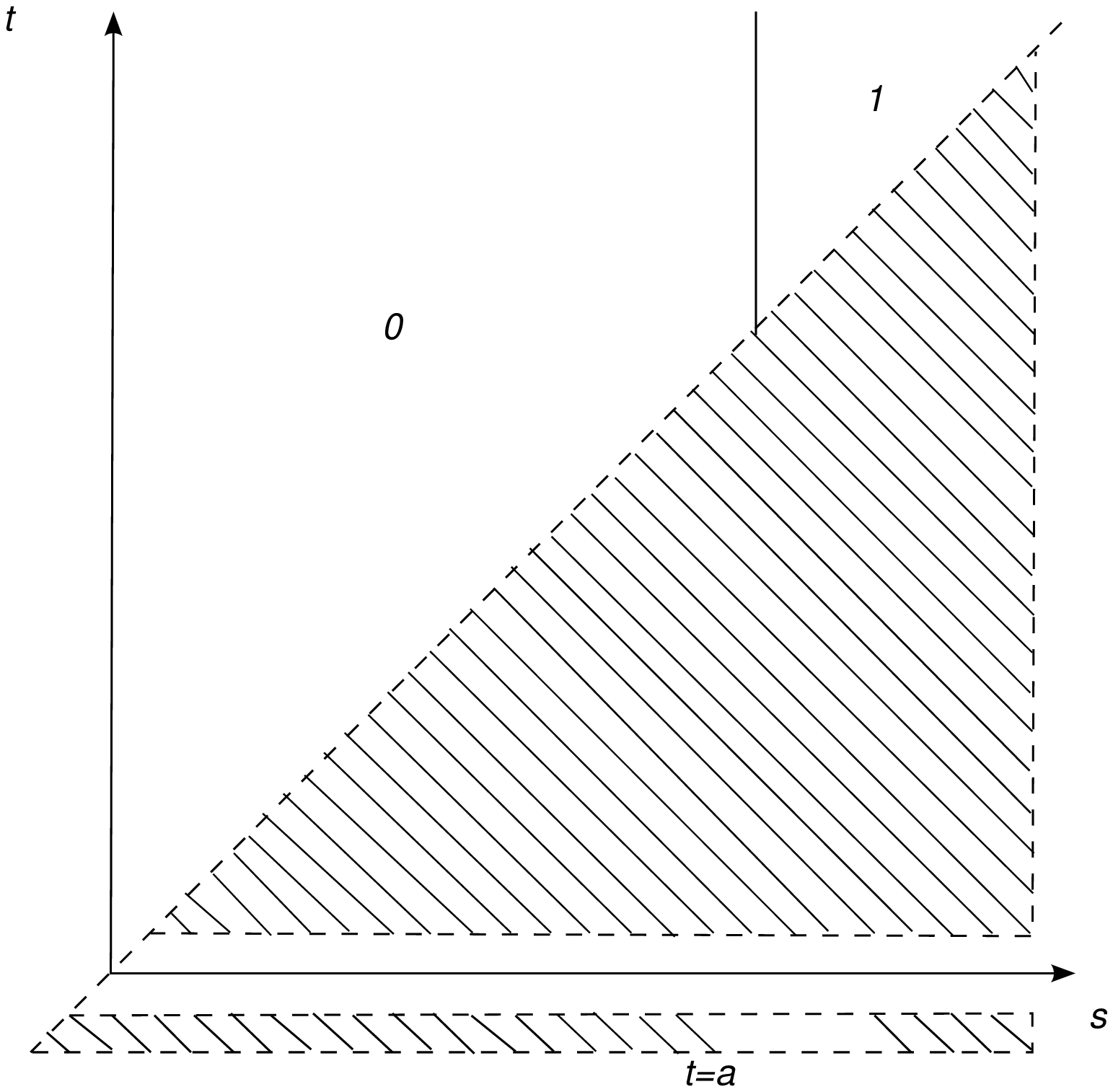}
\end{tabular}\\
$G_2^{s,t} (\mathcal S) = \left\{
\begin{array}{ll}
k,&  \sqrt{2}\le s<t\\
0, & \mbox{otherwise}
\end{array}
\right.$&
\begin{tabular}{c}
\psfrag{0}{$\scriptstyle{0}$} \psfrag{1}{$\scriptstyle{1}$}
\psfrag{2}{$\scriptstyle{2}$} \psfrag{s}{$\scriptstyle{s}$}
\psfrag{t}{$\scriptstyle{t}$}
\psfrag{t=a}{$\scriptstyle{\sqrt{2}}$}
\includegraphics[height=2cm]{Gst2SC.eps}
\end{tabular}\\
\end{tabular}
\end{array}\right\}$ &$ \Rightarrow D(\rho_{\C,2},\rho_{\mathcal S,2})\geq \frac{\sqrt{2}}{2} d(\rho'_{\C,2}\ \rho''_{\mathcal S,2}) =0$\\
\end{tabular}

In other words, multidimensional persistent homology, with respect
to $\vhr$ and $\psr$, is able to discriminate the cube and the
sphere, while the 1-dimensional one, with respect to $\f_1, \f_2$
and $\psi_1, \psi_2$, cannot do that. In fact, for either manifold
the lower level sets of the single components (i.e. 1-dimensional
measuring functions) are homeomorphic for all values: they are
topologically either circles, or annuli, or spheres.

It should be noted that the map $g''$ on
$\mathcal S$ reaches the homological critical value 1 at points,
at which it lacks of differentiability.

In the example above, $\vhr$ is not a Morse function (as would be
desirable, if not necessary), because
$g': \C \rightarrow \R$ has infinitely many critical points when
$\max\{|u|, |v|\} = 1$; moreover, the cubic surface itself is not
even $\mathcal C^1$. This problem can be solved by perturbing $\C$
so that it becomes smooth (e.g. a super-quadric \cite{JaLeSo00}).
In this case, the differences between homology modules of the cube
and of the super--quadric are only quantitative (i.e. the levels
of homological critical values are different from one another).

An even simpler example can be given on size pairs having the same
support. Let $X$ be the ellipse imbedded in $\R^3$ as
$\begin{cases} u^2+v^2=1\\v=w\end{cases}$ --- or parameterized as
$\begin{cases} u=\cos\theta\\
v=\sin\theta\\w=\sin\theta \end{cases}$. Let
$\f_1,\f_2,\psi_1,\psi_2:X\to\R$ be defined as $\f_1=u,\ \psi_1=v,
\ \f_2=\psi_2=w$ and $\vhr=(\f_1,\f_2),\ \psr=(\psi_1,\psi_2)$.
Then the persistent homology modules of $(X,\f_1)$, $(X,\psi_1)$,
$(X,\f_2)=(X,\psi_2)$ are identical, while the persistent homology
(in degree zero, so the size function) of $(X,\vhr)$ differs from
the one of $(X,\psr)$. Indeed, while the lower level sets of
$\psr$ are always either empty or connected, the lower level sets
$\vhr\le (\overline{u},\overline{w})$, with $0<\overline{u}<1,
\sqrt{1-\overline{u}^2}\le\overline{w}<1$ consist of two connected
components.

\section{Reduction of $i$-essential critical values}\label{essentiality}

The former example of the previous section suggests also some
other considerations on the cooperation of measuring functions.
We remind that the adjective ``essential'' is used here with the meaning
introduced in \cite{CaFePo01} and recalled in Section
\ref{1Dpersistent} after Definition \ref{homcrit}, so otherwise
than, e.g., in \cite{ChFr}.

A first remark is that, although the persistent homology on single components of $\fr$
cannot distinguish the two spaces, the persistent homology on $f_1$
restricted to lower level sets of $f_2$ can, as can be shown as follows.
Consider again the sphere $\mathcal S$. The value $1/\sqrt{2}$
(corresponding to the homological critical value 1 of $g''$) is not
critical for the maps $f_1, f_2$ on $\mathcal S$ itself, but it is
indeed critical for $f_2$ restricted to $f_1^{-1}(-\infty,1/\sqrt{2}]$.
 We believe that homological
critical values of the 1D reduction of multidimensional measuring
functions are always clues of such phenomena.

A further speculation on the use of cooperating measuring
functions --- from a completely different viewpoint than the one
developed in the previous sections --- is the following. A problem
in 1-dimensional persistent homology, as well as for the size
functor, is the computation of $i$-essential critical values for
$i>0$. A possibility is the use of several, independent measuring
functions for lowering $i$, i.e. the degree at which the passage
through the critical value causes a homology change. Lowering $i$
is important, since 0-essential critical values are easily
detected by graph-theoretical techniques \cite{DA00}. The
following example shows that a suitable choice of a second,
auxiliary measuring function may actually take 1-essential
critical values to 0-essential ones.

Let $\mathcal T$ be a torus of revolution around the $x$ axis, with
the innermost parallel circle of radius 2, the outermost of radius
3. On $\mathcal T$ define $(f_1,f_2)=(z,-z)$. Suppose we are interested
in the persistent homology of the size pair $(\mathcal T, f_1)$. Then
$(0,0,2)$ is a 1-essential critical point for $f_1$, i.e. it is a point
at which 1-degree homology changes. Of course, there are computational methods
(e.g. by the Euler--Poincar\'e characteristic) which enable us to detect it,
but they will probably be tailored to the particular dimension of the manifold
and to the particular homology degree.

The same point is 0-essential for
its restriction to $f_2^{-1}(-\infty,1]$, so it can be recovered by the standard
graph-theoretical techniques used in degree 0, i.e. for size functions. (The two functions need
not be so strictly related: $f_2$ could be replaced by Euclidean
distance from $(0,0,3)$ with the same effect). We conjecture that
--- at least whenever torsion is not involved --- one can
recursively take the $i$-essential values of a measuring function
to $(i-1)$-essential ones, down to (easily computable) 0-essential
critical values by means of other (auxiliary) measuring functions,
as in this example.

\section{Conclusions and future work}
The need of extending persistent homology to the multidimensional
case is a rather widespread belief, confirmed by simple examples.
The present research shows the possibility of reducing the
computation of persistent homology, with respect to
multidimensional measuring functions, to the 1-dimensional case,
following the line of thought of an analogous extension devised
for size functions in \cite{BiCeXX}. This reduction also yields a
stable distance for the rank invariants of size pairs.

In the next future, we plan to characterize the multidimensional
max-tame measuring functions in a way that the reduction to 1D
case makes the specific features of persistent homology modules
hold steady. It also would be our concern to give a rigorous
definition of \emph{multidimensional homological critical values}
of a max-tame function and to relate them to the homological
critical values of the maximum of its components.

Eventually, in relation to our conjecture about $i$-essentiality
(see Section \ref{essentiality}), we plan to build an algorithm to
recursively reduce $i$-essential critical points of a measuring
function to 0-essential ones.

\subsection*{Acknowledgements}
Work performed under the auspices of INdAM-GNSAGA, CIRAM, ARCES
and the University of Bologna, funds for selected research topics.


\begin{thebibliography}{100}


\bibitem{BiCeXX}
S. Biasotti, A. Cerri, P. Frosini, D. Giorgi and C. Landi, \emph{
Multidimensional size functions for shape comparison}, Journal of
Mathematical Imaging and Vision (in press).

\bibitem{CaFePo01}
F.~Cagliari, M.~Ferri and P.~Pozzi, \emph{Size functions from a
categorical viewpoint}, Acta Appl. Math. \textbf{67} (2001),
225-235.

\bibitem{CaZo06}
G. Carlsson and A. Zomorodian, \emph{The Theory of
Multidimensional Persistence}, Symposium on Computational Geometry,
June 6--8, 2007, Gyeongiu, South Korea (2007) 184--193.

\bibitem{ChFr}
C. Chen and D. Freedman, \emph{Quantifying homology classes},
25th Symp. on Theoretical Aspects of Computer Science, Bordeaux, France (2008), 169-180

\bibitem{CoEdHa05}
D.~Cohen-Steiner, H.~Edelsbrunner, J.~Harer, \emph{Stability of
persistence Diagrams}, Proc. 21st Sympos. Comput. Geom. (2005),
263--271.

\bibitem{DA00}
M. d'Amico, \emph{$\Delta^*$ reduction of size graphs as a new algorithm for computing size functions
of shapes} In: Proc. Internat. Conf. on Computer Vision, Pattern Recognition and Image
Processing, Feb. 27–Mar. 3, 2000, Atlantic City, vol. 2 (2000), 107–-110.

\bibitem{DAFrLa}
M. d'Amico, P. Frosini and C. Landi, \emph{Natural pseudo-distance
and optimal matching between reduced size functions}, Acta
Applicandae Mathematicae (to appear).

\bibitem{DeGh}
V. De Silva and R. Ghrist, \emph{Homological Sensor Networks},
Notices Amer. Math. Soc., \textbf{54}, no. 1 (2007), 10--17.

\bibitem{DeRGh}
V. De Silva and R. Ghrist, \emph{Coverage in sensor networks via
persistent homology}, Alg. and Geom. Topology, 7, (2007) 339--358.

\bibitem{EdLeZo00}
H.~Edelsbrunner, D.~Letscher and A.~Zomorodian, \emph{Topological
Persistence and Simplification}, Proc. 41st Ann. IEEE Sympos.
Found Comput. Sci. (2000), 454--463.

\bibitem{EdLeZo02}
H.~Edelsbrunner, D.~Letscher and A.~Zomorodian, \emph{Topological
Persistence and Simplification}, Discrete Comput. Geom.
\textbf{28} (2002), 511--533.

\bibitem{Fr91}
P.~Frosini, \emph{Measuring shapes by size functions}, Proc. of
SPIE, Intelligent Robots and Computer Vision X: Algorithms and
Techniques, Boston, MA 1607 (1991), 122--133.

\bibitem{FrLa01}
P.~Frosini and C.~Landi, \emph{Size functions and formal series},
Appl. Algebra Eng. Commun. Computing \textbf{12} (2001), 327--349.

\bibitem{FrMu99}
P.~Frosini and M.~Mulazzani, \emph{Size homotopy groups for
computation of natural size distances}, Bull. Belg. Math. Soc.
\textbf{6} (1999), 455--464.

\bibitem{Gh08}
R. Ghrist, \emph{Barcodes: The persistent topology of data},
Bull. Amer. Math. Soc. \textbf{45} (2008), 61--75.

\bibitem{JaLeSo00}
A. Jaklic, A. Leonardis and F. Solina, \emph{Segmentation and
Recovery of Superquadrics}, Computational imaging and vision,
Kluwer, Dordrecht, \textbf{20} (2000), ISBN 0-7923-6601-8.

\bibitem{KaMiMr04}
T. Kaczynski, K. Mischaikow  and M. Mrozek, {Computational
Homology}, Applied Mathematical Sciences \textbf{157},
Springer-Verlag, New York (2004).

\bibitem{LaFr97}
C.~Landi and P.~Frosini, \emph{New pseudodistances for the size function space},
Proc. SPIE Vol. 3168, Vision Geometry VI, Robert A. Melter, Angela Y.
Wu, Longin J. Latecki (eds.) (1997), 52--60.

\bibitem{UrVe97}
C.~Uras and A.~Verri, \emph{Computing size functions from edge
maps}, Internat. J. Comput. Vision \textbf{23} (1997), no.~2,
169--183.

\bibitem{VeUr96}
A.~Verri and C.~Uras, \emph{Metric-topological approach to shape
representation and recognition}, Image Vision Comput. \textbf{14}
(1996), 189--207.

\bibitem{VeUrFrFe93}
A.~Verri, C.~Uras, P.~Frosini and M.~Ferri, \emph{On the use of
size functions for shape analysis}, Biol. Cybern. \textbf{70}
(1993), 99--107.

\end{thebibliography}
\end{document}